\documentclass[12pt]{amsart}
\usepackage{amsthm,amsmath,amssymb}
\usepackage{fullpage}
\usepackage{hyperref}
\usepackage{cleveref}
\usepackage{lmodern}
\usepackage{xcolor}
\usepackage{setspace}
\usepackage{todonotes}
\hypersetup{
  colorlinks,
  linkcolor={blue!80!black},
  citecolor={blue!80!black},
  urlcolor={blue!100!black}
}

\theoremstyle{plain}
\newtheorem{lemma}{Lemma}[section]
\newtheorem{theorem}{Theorem}[section]
\newtheorem*{theorem*}{Theorem}
\newtheorem{corollary}{Corollary}[section]
\newtheorem{proposition}{Proposition}[section]

\theoremstyle{definition}
\newtheorem{definition}{Definition}[section]

\newtheorem{example}{Example}[section]
\newtheorem*{question}{Question}

\theoremstyle{remark}
\newtheorem*{remark}{Remark}

\crefname{theorem}{theorem}{theorems}
\Crefname{theorem}{Theorem}{Theorems}
\crefname{lemma}{lemma}{lemmas}
\Crefname{lemma}{Lemma}{Lemmas}
\Crefname{proposition}{Proposition}{Propositions}
\Crefname{corollary}{Corollary}{Corollaries}

\newcommand{\abs}[1]{\left\lvert #1 \right\rvert}

\newcommand{\BF}{\mathbb{F}}
\newcommand{\BZ}{\mathbb{Z}}
\newcommand{\FT}{\mathbb{F}_2}

\newcommand{\mf}{\mathfrak}
\newcommand{\mc}{\mathcal}
\DeclareMathOperator{\Ann}{Ann}
\DeclareMathOperator{\Spec}{Spec}
\DeclareMathOperator{\Stab}{Stab}

\newcommand{\intersect}{\cap}
\begin{document}
\title{Davenport Constant for Commutative Rings}
\author{Calvin Deng}
\begin{abstract}
  The Davenport constant is one measure for how ``large'' a finite abelian group
  is.  In particular, the Davenport constant of an abelian group is the smallest
  $k$ such that any sequence of length $k$ is reducible.  This definition
  extends naturally to commutative semigroups, and has been studied in certain
  finite commutative rings. In this paper, we give an exact formula for the
  Davenport constant of a general commutative ring in terms of its unit group.
\end{abstract}
\maketitle
\section{Introduction}
The Davenport constant is an important concept in additive number theory. In 
particular, it measures the largest zero-free sequence of an abelian group.

The Davenport constant was introduced by Davenport in 1966 \cite{Davenport}, but
was actually studied prior to that in 1963 by Rogers \cite{Rogers}.
The definition was first extended by Geroldinger and Schneider to abelian semigroups by \cite{GHK2006} as follows:
\begin{definition}
  For an additive abelian semigroup $S$, let $d(S)$ denote the smallest $d\in
  \mathbb{N}_0\cup \{\infty\}$ with the following property:

  For any $m\in \mathbb{N}$ and $s_1, \dots, s_m\in S$ there exists a subset
  $J\subseteq [1, m]$ such that $\abs{J}\le d$ and
  \[
    \sum_{j = 1}^m s_j = \sum_{j\in J} s_j.
  \]
\end{definition}
In addition, \cite{GHK2006} showed the following:
\begin{proposition}
  \label{infprop}
  If $\abs{S} < \infty$, then $d(S) < \infty$.
\end{proposition}
Wang and Guo \cite{WangGuo2008} then gave the definition of the large 
Davenport constant in terms of reducible and irreducible sequences, as follows:
\begin{definition}
  Let $S$ be a commutative semigroup (not necessarily finite). Let $A$ be a
  sequence of elements in $S$. We say that $A$ is \textit{reducible} if there exists a
  proper subsequence $B\subsetneq A$ such that the sum of the elements in $A$ is
  equal to the sum of the elements in $B$.  Otherwise, we say that $A$ is
  \textit{irreducible}.
\end{definition}
\begin{definition}
  Let $S$ be a finite commutative semigroup. Define the \textit{Davenport constant}
  $D(S)$ of $S$ as the smallest $d\in \mathbb{N} \cup \{\infty\}$ such that
  every sequence $S$ of $d$ elements in $S$ is reducible.
\end{definition}
\begin{remark}
  $D$ and $d$ are related by the equation $D(S) = d(S) + 1$.
\end{remark}
Note that if $S$ is an abelian group, being irreducible is equivalent to being
zero-sum free, so the definition of the Davenport constant here is equivalent to 
the classical definition of the Davenport constant for abelian groups.

In all following sections, unless otherwise noted:
\begin{itemize}
  \item All semigroups are unital and commutative.
Furthermore, we will use multiplication notation for semigroups, as opposed to
the additive convention used in \cite{GHK2006,WangGuo2008,Wang2015124,QuWangZhang}.
  \item Similarly, rings are unital and commutative.
  \item Sets represented by the capital letter $S$ are semigroups.
  \item Sets represented by the capital letter $T$ are ideals in a
    semigroup.
  \item Sets represented by the capital letters $A, B$ are
    sequences of elements in a semigroup. In addition, $\pi(A)$ denotes the
    product of the elements in $A$.
  \item Sets represented by the capital letter $R$ are commutative rings.
  \item By abuse of notation, when we write $D(R)$, we actually mean $D(S_R)$,
      where $S_R$ is the semigroup of $R$ {\it under multiplication}.
  \item $C_n$ denotes
    the cyclic group of order $n$; $\BZ / n\BZ$ denotes the ring with additive group $C_n$.
\end{itemize}
\section{Previous Results}
So far, this general setting has been studied extensively in the cases where the 
semigroup is the semigroup of a (commutative) ring under multiplication.
The general idea in this class of problems is to show that the Davenport
constant of the semigroup of a ring $R$ under multiplication is ``close'' to
the Davenport constant of its group of units $U(R)$. The reason this has been done
is that very little is known about the precise value of $D(G)$ when $G$ is an
abelian group of high rank, and often the unit group of these commutative rings
will be a group of high rank. However, as we will show, if there is a general theorem
for the Davenport constant of an abelian group, then by \Cref{thm1}, we will also have
a general theorem for the Davenport constant of an arbitrary finite commutative ring.
\begin{remark}
Since $U(R)$ 
is a sub-semigroup of $R$, we clearly have $D(R)\ge D(U(R))$, so we would like to 
say something about the difference $D(R) - D(U(R))$.
\end{remark}
In their paper, Wang and Guo showed the following result:
\footnote{\cite{WangGuo2008} actually erroneously claimed that for general
$n_i$, $D(R) = D(U(R)) + \#\{1\le i\le r: 2\| n_i\}$. The author has
corresponded with the authors of \cite{WangGuo2008} on this matter, and they
offer the result given above.}
\begin{theorem}[Wang, Guo, 2008, \cite{WangGuo2008}]
  \label{wangguo}
  If $R = \BZ/n_1\BZ \times \dots \times \BZ/n_r\BZ$, where each of the 
  $n_1, \dots, n_r$ are odd. Then 
  $D(R) = D(U(R))$.
\end{theorem}
Later, Wang showed the following theorem:
\begin{theorem}[Wang, 2015, \cite{Wang2015124}]
  \label{wang}
  Suppose $q > 2$ is a prime power.
  If $R \neq \BF_q[x]$ is a quotient ring of $\BF_q[x]$, then $D(R) = D(U(R))$.
\end{theorem}
However, Wang left the case $q = 2$ open and gave an instance for which 
$D(R)\neq D(U(R))$. Later, Zhang, Wang, and Qu gave the following bound
when $q = 2$:
\begin{theorem}[Zhang, Wang, Qu, 2015, \cite{QuWangZhang}]
  \label{qwz}
  If $f\in \FT[x]$ is nonconstant and $R =  \FT[x] / (f)$,
  then $D(U(R))\le D(R) \le D(U(R)) + \delta_f$, where $\delta_f = \deg[\gcd(f,
  x(x + 1))]$.
\end{theorem}

\section{Summary of New Results}
First, we will show that the converse of \Cref{infprop} holds for rings:
\begin{theorem}
    \label{thm4}
  If $R$ is a commutative ring and $D(R) < \infty$, then $\abs{R} < \infty$.
\end{theorem}
However, the main result of this paper will be to relate $D(R)$ and $D(U(R))$ for arbitrary finite rings $R$.

As we will see, the reason why $D(R) = D(U(R))$ does not hold in general is closely related to the presence of certain index 2 ideals in $R$, which we can see immediately in the following exact formula for $D(R)$ in terms of the Davenport constants of certain subgroups of $U(R)$:
\begin{theorem}
  \label{thm1}
  Suppose $R$ is of the form
  \[
    (\BZ / 2\BZ)^{k_1} \times (\BZ / 4\BZ)^{k_2} \times (\BZ / 8\BZ)^{k_3} \times (\FT[x] / (x^2))^{k_4
    } \times R',
  \]
  where $R'$ is a product of local rings not isomorphic to $\BZ / 2\BZ, \BZ / 4\BZ, \BZ / 8\BZ, \FT[x]/(x^2)$.
  Then
  \[
    D(R) = \max_{0\le a\le k_2 + k_4, 0\le b\le k_3}\bigg[D(U(R') \times {C_2}^{k_2 + k_4 + 2k_3 - a -
    2b}) + k_1 + 2a + 3b \bigg].
  \]
\end{theorem}
\begin{remark}
  Note that since $R$ is finite, $R$ is Artinian, so there is a unique decomposition of $R$ as a product as a product of local rings (\cite{atiyahmacdonald}, \S 8). Thus the quantities $k_1, k_2, k_3, k_4$ are well-defined as functions of $R$.
\end{remark}
As a corollary, we get the following bound on $D(R) - D(U(R))$:
\begin{corollary}
    \label{thm2}
  Suppose $R$ is of the form
  \[
    (\BZ / 2\BZ)^{k_1} \times (\BZ / 4\BZ)^{k_2} \times (\BZ / 8\BZ)^{k_3} \times (\FT[x] / (x^2))^{k_4
    } \times R',
  \]
  where $R'$ is a product of local rings not isomorphic to $\BZ / 2\BZ, \BZ / 4\BZ, \BZ / 8\BZ, \FT[x]/(x^2)$.
  Then
  \[
    D(U(R)) + k_1 \le D(R) \le D(U(R)) + k_1 + k_2 + k_3 + k_4.
  \]
  In addition, equality holds on the right if $U(R)$ is a power of 2 or if
  $k_i = 0$ for $i = 2, 3, 4$.

  \begin{proof}
    For the left hand side, note that we can pick $a = b = 0$ to get
    \[
      D(R)\ge D(U(R') \times {C_2}^{k_2 + k_4 + 2k_3}) + k_1 = D(U(R)) + k_1
    \]
    For the right hand side, we use the well-known facts $D(G\times H) \ge D(G) + D(H) - 1$, and $D(C_2) = 2$ to get
    \begin{align*}
      D(R) &= \max_{0\le a\le k_2, 0\le b\le k_3}\bigg[D(U(R') \times {C_2}^{k_2 + k_4 + 2k_3 - a -
      2b}) + k_1 + 2a + 3b \bigg] \\ &= \max_{0\le a\le k_2, 0\le b\le
      k_3}\bigg[D(U(R') \times {C_2}^{k_2 + k_3 + 2k_3 - a - 2b}) + (a + 2b) (D(C_2) - 1) + k_1 + a
      + b \bigg] \\
      &\le \max_{0\le a\le k_2, 0\le b\le k_3}\bigg[D(U(R') \times
        {C_2}^{k_2 + k_4 + 2k_3}) + k_1 + a + b \bigg] \\
      &= \max_{0\le a\le k_2, 0\le b\le k_3}\bigg[D(U(R)) + k_1 + a + b \bigg] \\
      &= D(U(R)) + k_1 + k_2 + k_3 + k_4.
    \end{align*}
    When $k_2 + k_3 + k_4 = 0$, equality clearly holds as the left hand side and the right hand side are the same.
    On the other hand, when $\abs{U(R)}$ is a power of 2,
    then $U(R')$ is also a 2-group. However, by \cite{Olson1969}, if $G$ and $H$ are
    2-groups, then $D(G\times H) = D(G) + D(H) - 1$. Thus
    \begin{align*}
      D(R) &\ge D(U(R')) + k_1 + 2k_2 + 3k_3 + 2k_4 \\
      &= D(U(R')) + (k_2 + k_4 + 2k_3) (D(C_2) - 1) + k_1 + k_2 + k_3 + k_4 \\
      &= D(U(R')\times {C_2}^{k_2 + k_4 + 2k_3}) + k_1 + k_2 + k_3 + k_4 \\
      &= D(R) + k_1 + k_2 + k_3 + k_4
      \ge D(R),
    \end{align*}
    so equality holds in this case as well.
  \end{proof}
\end{corollary}
We also have the following more concise bound on $D(R) - D(U(R))$ that does
not depend on writing down a local ring product decomposition.
\begin{corollary}
    \label{thm3}
  Suppose $R$ is a finite commutative ring, and let $n_2(R)$ be the number of
  index two (prime, maximal) ideals of $R$. Then $D(R)\le D(U(R)) + n_2(R)$.
  \begin{proof}
    Since $R$ is finite (and thus Artinian), $R$ can be expressed as a finite
    product of local rings. Thus in the setting of \Cref{thm2}, it suffices to
    show that $k_1 + k_2 + k_3 + k_4\le n_2(R)$. Since $R$ is Artinian, we have that $R\simeq
    \prod_{i = 1}^n R / \mf{q}_i$, where $\{\sqrt{\mf{q}_i}\}_{i = 1}^n =
    \Spec(R)$. In particular, the number of $R / \mf{q}_i$ that can be
    isomorphic to $\BZ / 2\BZ, \BZ/4\BZ, \BZ/8\BZ$ or $\FT[x] / (x^2)$ is at
    most the number of index 2 ideals in $\Spec(R)$, which is precisely
    $n_2(R)$.
  \end{proof}
\end{corollary}
Using \Cref{thm1} and \Cref{thm2,thm3}, we will give generalizations of each of the results in the previous section.

The main line of attack to prove \Cref{thm1} will be to reduce the problem to the case of finite local rings and then talk about what happens when we ``glue'' local rings together via product. However, we will first discuss the gluing mechanism for a more general class of semigroups (``almost unit-stabilized''); this will require the notion of \textit{unit-stabilized pairs} and a \textit{relative Davenport constant}. Afterwards, we will show that finite local rings are almost unit-stabilized, and that more structure holds for all finite local rings other than $\BZ / 2\BZ, \BZ / 4 \BZ, \BZ / 8\BZ$, and $\FT[x] / (x^2)$. 

\section{Davenport constant for semigroups}
Let $S$ be a semigroup and $U(S)$ denote its group of units.
\begin{definition}
  We say that a sequence $A$ with elements in
  $S$ is \textit{reducible} if there exists a proper subsequence $A'\subsetneq A$ such that
  $\pi(A') = \pi(A)$. Otherwise, $A$ is \textit{irreducible}.
\end{definition}
\begin{definition}
  The Davenport constant $D(S)$ is the minimum positive integer 
  $k$ such that any sequence $A$ of size $k$ is reducible.
\end{definition}
\begin{definition}
  For $S$ a semigroup, let $\Delta(S) = D(S) - D(U(S))$.
\end{definition}
\begin{remark}
  Clearly $\Delta(S) \ge 0$.
\end{remark}
\begin{lemma}
  \label{prodlemma}
  Suppose $S$ and $S'$ are semigroups. Then
  \[
    D(S\times S')\ge D(S) + D(S') - 1.
  \]
  \begin{proof}
    Let $\ell = D(S) - 1$ and $\ell' = D(S') - 1$. Then there exists an irreducible sequence $s_1, \dots, s_\ell$ in $S$ of length $\ell$ and an irreducible sequence $s_1',\dots s_{\ell'}'$ in $S'$ of length $\ell'$. Then the sequence
    \[
      (s_1, 1_{S'}), \dots, (s_\ell, 1_{S'}),
      (1_S, s_1'), \dots, (1_S, s_{\ell'}')
    \]
    is an irreducible sequence in $S\times S'$ of length $\ell + \ell' = D(S) + D(S') - 2$, which means that $D(S\times S') \ge D(S) + D(S') - 1$.
  \end{proof}
\end{lemma}
\section{Proof of \texorpdfstring{\Cref{thm4}}{Theorem \ref{thm4}}}
We will first handle the case where $R$ is infinite. In particular, we
will show that if $\abs{R} = \infty$, then $D(R) = \infty$.
\begin{remark}
  The corresponding theorem for semigroups clearly does not hold: let $S$ be the
  semigroup with underlying set $\BZ^{\ge 0}$ and operation $a * b = \max(a,
  b)$. Then $S$ is infinite but $D(S) = 2$.
\end{remark}
\begin{lemma}
  \label{inflemma}
  Suppose $G$ is an infinite abelian group. Then $D(G) = \infty$.
  \begin{proof}
    Let $g_1, \dots, g_n$ be an irreducible sequence. Then $g_1, \dots, g_n, g$
    is not irreducible if and only if $g^{-1}$ is not of the form $g_{i_1} \cdots g_{i_k}$.
    However, there are at most $2^n$ such $g$, which means that we can always
    extend an irreducible sequence to get a longer irreducible sequence, which
    means that $D(G) = \infty$.
  \end{proof}
\end{lemma}
\begin{lemma}
  Suppose $R$ is a ring and $I\subseteq R$ is an ideal. Then $D(R)\ge D(R / I)$.
  \begin{proof}
    This follows from the fact that any irreducible sequence in $R / I$ can
    be lifted to an irreducible sequence in $R$.
  \end{proof}
\end{lemma}
We are now ready to give the proof of \Cref{thm4}:
\begin{theorem*}
  If $R$ is a commutative ring and $D(R) < \infty$, then $\abs{R} < \infty$.
  \begin{proof}
    Suppose $D(R) = k < \infty$.
    \begin{itemize}
      \item $\abs{R / \mf{m}}$ is finite for every maximal ideal $\mf{m}
        \subseteq R$. Otherwise, $U(R / \mf{m})$ is an infinite abelian
        group, so $D(R / \mf{m}) = \infty$, which means that $D(R) =
        \infty$ as well.
      \item $R$ has at most $k - 1$ maximal ideals. Suppose otherwise, that
        $R$ has $k$ maximal ideals $\mf{m}_1, \dots, \mf{m}_k$. Then by
        the Chinese Remainder theorem, there exists $a_1, \dots, a_k\in
        R$ such that $a_i \in \mf{m}_i$ and $a_i - 1 \in \mf{m}_j$ for
        $i\neq j$, and $a_1, \dots, a_k$ would be an irreducible
        sequence of length $k$, contradiction.
      \item Let $J(R)$ denote the Jacobson radical of $R$, the
        intersection of all of the maximal ideals of $R$. Since $R$
        has a finite number of maximal ideals,
        by the Chinese Remainder Theorem, $R / J(R)\simeq \prod_{i}
        R / \mf{m}_i$, which means that $R / J(R)$ is finite.
      \item On the other hand, $\abs{U(R)} < \infty$, for if $\abs{U(R)} = \infty$, then $D(U(R)) = \infty$, so $D(R) =
        \infty$ as well.
      \item Finally, $1 + J(R)\subseteq U(R)$, which means that
        $\abs{J(R)} < \infty$ as well. Thus $\abs{R} = \abs{R / J(R)}
        \cdot \abs{J(R)} < \infty$.
    \end{itemize}
  \end{proof}
\end{theorem*}
\section{Relative Davenport constant}
In order to prove \Cref{thm1} (and its subsequent corollaries), we need to introduce the notion of the \textit{relative Davenport constant}, which measures the longest irreducible sequence in $S$ with product lying in an ideal $T\subseteq S$.
\begin{definition}
    A set $T\subseteq S$ is an ideal if $S\cdot T\subseteq T$. ($T = \emptyset$
    is allowed here.)
\end{definition}
Suppose $S_1, S_2$ are semigroups and $T_1\subseteq S_1, T_2\subseteq S_2$ are
ideals. Then $T_1\times T_2\subseteq S_1\times S_2$ is also an ideal.
\begin{definition}
  Given a semigroup $S$ and an ideal $T\subseteq S$, $A$ is a \textit{$T$-sequence}
  if $\pi(A)\in T$.
\end{definition}
\begin{definition}
  Given a semigroup $S$ and a (possibly empty) ideal $T\subseteq S$, define 
  the \textit{relative Davenport constant} $D(S, T)$ to be the minimum positive integer 
  $k$ such that any $T$-sequence $A$ of size $k$ is reducible.
\end{definition}
We can again define $d(S, T)$ as the maximum length of an irreducible $T$-sequence, and we once again have the relation $D(S, T) = d(S, T) + 1$.
\begin{remark}
  If $T = \emptyset$ is the empty ideal, then $d(S, T)$ is tautologically 0 and so $D(S, T)$ is 1.
\end{remark}
Before we proceed, we first need the following helpful lemma.
\begin{lemma}
  \label{rellemma}
  Let $S$ be a semigroup and $T$ be an ideal. Suppose $s\in S$ and $A$ is a 
  sequence in $S$ of length $D(S, T)$ such that $s\cdot \pi(A) \in T$. Then 
  there exists a proper subsequence $A'\subsetneq A$ such that $s\pi(A') = 
  s\cdot \pi(A)$.
  \begin{proof}
    If $\pi(A)\in T$, then we are done by the definition of $D(S, T)$. 
    Otherwise, consider the sequence $B = A\cup \{s\}$. Since $\abs{B} = 
    D(S, T) + 1 > D(S, T)$, there exists a proper subsequence $B'\subsetneq 
    B$ such that $\pi(B') = \pi(B) = s\cdot \pi(A)$. Note that if $s\notin B'$, 
    then $B'\subseteq A$, which means that $\pi(B')\notin T$ (as 
    $\pi(A)\notin T$ by definition), whereas $\pi(B) \in T$, contradiction. 
    Thus $a\in B'$, so if we let $A' = B' \backslash \{s\}$, then $s\cdot \pi(A') 
    = \pi(B') = \pi(B) = s\cdot \pi(A)$, as desired.
  \end{proof}
\end{lemma}
\begin{lemma}
  \label{grplem}
  Suppose $G$ is an abelian group and $H$ is a subgroup of $G$, and suppose 
  $S$ is a semigroup.
  Then
  \[
    D(G \times S, G \times T) \ge D((G / H) \times S, (G / H) \times T) + 
    D(H) - 1.
  \]
  \begin{proof}
    Let $\ell_1 = D(G / H\times S, (G / H)\times T) - 1, \ell_2 = D(H) - 1$. 
    Then there exists an irreducible ($(G/ H)\times T$)-sequence $(g_1, 
    s_1), \dots (g_{\ell_1}, s_{\ell_1})$ in $(G / H) \times S$ of length 
    $\ell_1$, and an irreducible sequence $h_1, \ldots, h_{\ell_2}$ in $H$ 
    of length $\ell_2$. In this case, if we let $g_i'$ be any lift of $g_i$ 
    from $G / H$ to $G$, then the sequence
    \[
      (g_1', s_1), \ldots, (g_{\ell_1}', s_{\ell_1}), (h_1, 1_S), \ldots, 
      (h_{\ell_2}, 1_S)
    \]
    is an irreducible $G\times T$-sequence in $G\times S$ of length
    \[
      \ell_1 + \ell_2 = D((G / H)\times S, (G / H)\times T) + D(H) - 2.
    \]
  \end{proof}
\end{lemma}
\section{Unit-Stabilized Pairs}
\begin{definition}
  An ordered pair $(S,T)$ of a semigroup $S$ and an ideal $T\subseteq S$ is a 
  \textit{unit-stabilized pair} if for all $a, b\in S$ such that $ab\notin T$ 
  and $\Stab_{U(S)}(a) = \Stab_{U(S)}(ab)$, then $ab = au$ for some $u \in U(S)$.
\end{definition}
\begin{definition}
  We say that $S$ is a \textit{unit-stabilized semigroup} if $(S, \emptyset)$ is a unit-stabilized pair.
\end{definition}
\begin{example}
  If $R$ is a finite field other than $\FT$, then [the semigroup of] $R$ [under
  multiplication] is unit-stabilized. $(\FT, 0)$ is a unit-stabilized pair.
\end{example}
It turns out most local rings $R$ are unit-stabilized.
Unit-stabilized rings are nice because they satisfy $\Delta(R) = 0$ and are closed
under the operation of taking products. However, if $R$ has residue field $\FT$, then 
it is in fact not unit-stabilized: $U(R) = 1 + \mf{m}$, where $\mf{m}$ is the maximal
ideal of $R$. Thus if $x\in R$ satisfies $\Ann_R(x) = \mf{m}$, then $\Stab_{U(R)}(x) = \Stab_{U(R)}(0) = U(R)$, but clearly $0\neq ux$ for any unit $u$.

As we can see, the unit-stabilized behavior breaks down around 0. However, 
even when $R$ is not unit-stabilized, it turns out that $(R, 0)$ will be a unit-stabilized pair.
We have the following formula that helps for working with unit-stabilized pairs:
\begin{theorem}
  \label{relthm}
  Suppose $(S, T)$ is a unit-stabilized pair and $S'$ is a semigroup and $T'$ 
  is an ideal in $S'$.
  Then
  \begin{equation*}
    D(S\times S', S\times T') = \max[D(U(S)\times S', U(S)\times T'), 
    D(S\times S', T \times T')].
  \end{equation*}
  \begin{proof}
    First off, it is clear that $D(S\times S', S \times T')\ge 
    D(U(S)\times S', U(S)\times T')$ and $D(S\times S', S\times T')\ge 
    D(S\times S', T\times T')$, so it suffices to show the inequality
    in the other direction.

    Let
    \[
      \ell = \max[D(U(S)\times S', U(S)\times T'), D(S\times S', T \times T')].
    \]
    Suppose we have a sequence $A = \{a_1,\ldots,a_\ell\}$ with $\pi(A)\in 
    S\times T'$, and suppose for the sake of contradiction that $A$ is 
    irreducible.

    Let $P_{S}$ and $P_{S'}$ denote the projection maps from $S\times S'$ to 
    $S$ and $S'$, respectively.
    If $P_S(\pi(A))\in T$, then $\pi(A)\in T\times T'$, and by the 
    definition of $D(S\times S', T \times T')$ we automatically have that 
    $A$ is not irreducible.
    Thus $P_S(\pi(A))\notin T$.

    Let $B$ be a minimal subsequence of $A$ such that 
    \[
      \Stab_{U(S)}(P_S(\pi(B))) = \Stab_{U(S)}(P_S(\pi(A))),
    \]
    and without loss of generality suppose that $B = \{a_1, \ldots, a_k\}$.  Let
    $p_0 = (1_{S}, 1_S)$, and for $1\le i\le \ell$, $p_i = a_ip_{i - 1}$ and let
    $K_i = \Stab_{U(S)}(P_S(p_i))$.  Then we have
    \[
      0 = K_0\subseteq K_1\subseteq \dots \subseteq K_\ell.
    \]
    In addition, for $1\le i\le k$, $K_{i-1}\neq K_i$, so by \Cref{grplem},
    \begin{align*}
      D(K_i)\ge D(K_{i - 1}) + D(K_i / K_{i - 1}) - 1 \ge D(K_{i - 1}) + 1,
    \end{align*}
    which means that $D(K_k)\ge k + 1$. Applying \Cref{grplem} again, we have, 
    \begin{align*}
      D((U(S) / K_k) \times S', (U(S) / K_k)\times T') \le
      D(U(S) \times S', U(S) \times T') - D(K_k) + 1
      \le \ell - k.
    \end{align*}
    In addition, for $k < i \le \ell$,
    \[
      \Stab(P_S(p_k)) \subseteq \Stab(P_S(a_ip_k)) 
    \subseteq \Stab(P_S(p_\ell)) = \Stab(P_S(p_k)),
    \]
    so the inclusions are equalities.
    Since $p_\ell = P_S(\pi(A))\notin T$ and $(S, T)$ is unit-stabilized,
    for each $k < i \le \ell$ there exists a $u_i\in U(S)$ such that
    $u_iP_S(p_k) = a_iP_S(p_k)$ for $k < i \le \ell$. Let $t_i = 
    (u_i, P_{S'}(a_i))$. We have $p_kt_i = p_ka_i$ for $k\le i\le \ell$.
    Since $t_{k+1}, \ldots, t_{\ell}$ is a sequence of length \[
      \ell - k\ge D((U(S) / K_k) \times S', (U(S) / K_k)\times T')
    \]
    and $(1, \pi_{S'}(p_k))\cdot t_{k+1}\cdots t_{\ell} \in U(S)\times T'$, by \Cref{rellemma},
    there exists a 
    (possibly empty) proper subsequence $i_1, \ldots, i_m$ such that $
    (1, \pi_{S'}(p_k))\cdot t_{i_1}\cdots t_{i_m}\cdot (u, 1_{S'}) = (1, \pi_{S'}(p_k))\cdot t_{k+1}\cdots t_{\ell}$ for 
    some $u\in K_k$.
    Then since $u\in \Stab(P_S(p_k))$, $(u, 1) \cdot p_k = p_k$, so
    \begin{align*}
      a_1a_2\cdots a_k \cdot a_{i_1}\cdots a_{i_m} &= p_k \cdot a_{i_1}\cdots a_{i_m} \\
      &= p_k \cdot t_{i_1}\cdots t_{i_m}\\
      &= (1, P_{S'}(p_k))\cdot (P_S(p_k), 1) \cdot t_{i_1}\cdots t_{i_m} \\
      &= (1, P_{S'}(p_k))\cdot (P_S(p_k), 1) \cdot (u, 1) \cdot t_{i_1}\cdots t_{i_m} \\
      &= (1, P_{S'}(p_k))\cdot (P_S(p_k), 1) \cdot t_{k + 1}\cdots t_{\ell} \\
      &= (1, P_{S'}(p_k))\cdot (P_S(p_k), 1) \cdot a_{k + 1}\cdots a_{\ell} \\
      &= p_k \cdot a_{k+1}\cdot a_\ell \\
      &= a_1\cdots a_\ell = \pi(A).
    \end{align*}
    This contradicts the assumption that $A$ was irreducible.
    Thus $D(S\times S', T \times T')\le \max[D(U(S)\times S', U(S)\times 
    T'), D(S\times S', T \times T')]$, as desired.
  \end{proof}
\end{theorem}
\begin{corollary}
  Suppose $(S, T)$ is a unit-stabilized pair. Then $D(S) = \max[D(U(S)), D(S, 
  T)]$.
  \begin{proof}
    Apply \Cref{relthm} to the case where $S'$ is the trivial group 
    and $T' = S'$.
  \end{proof}
\end{corollary}
\begin{corollary}
  If $S$ is a unit-stabilized semigroup, then $D(S) = D(U(S))$.
\end{corollary}

\section{Almost Unit-Stabilized semigroups}
\begin{lemma}
  \label{zerolemma}
  If $S$ has a 0 element, then for all semigroups $S'$,
  \[
    D(S \times S', \{0\}\times T') = D(S', T') + D(S, \{0\}) - 1.
  \]
  \begin{proof}
    First we will show that $D(S\times S', \{0\}\times T)\le D(S', T') + 
    D(S, \{0\}) - 1$.

    Let $\ell = D(S', T') + D(S, \{0\}) - 1$, and suppose we have an 
    irreducible $\{0\}\times T$-sequence $(s_1, s_1'), \ldots, (s_\ell, 
    s_\ell')$. Consider the smallest subset such that the second coordinate 
    of the product is equal to 0. By the definition of $D(S, \{0\})$,
    this subset has size at most $k 
    = D(S, \{0\}) - 1$. Without loss of 
    generality suppose $s_1\cdots s_k = 0$, and let $p = s_1'\cdots s_k'$. Then the sequence $s_{k+1}',\dots, 
    s_\ell'$ has length $D(S', T')$, so by \Cref{rellemma} there is some 
    proper subsequence such that $ps_{i_1}'\cdots s_{i_r}' = ps_{k + 
    1}'\cdots s_{\ell}'$. Since the first coordinate of both sides must be 
    zero, we have found a proper subsequence with the same product, 
    contradiction.

    Now to show the reverse inequality, it suffices to construct an 
    irreducible ($\{0\}\times T'$)-sequence of length $D(S', T') + D(S, 
    \{0\}) - 2$. Let $\ell_1 = D(S', T') - 1$ and $\ell_2 = D(S, \{0\}) - 
    1$. By definition, there exists an irreducible $T'$-sequence $s_1', 
    s_2', \ldots, s_{\ell_1}'\in S'$ and an irreducible $\{0\}$-sequence 
    $s_1, s_2, \ldots, s_{\ell_2}$. Then $(1_{S}, 
    s_1'), \dots, (1_{S}, s_{\ell_1}'), (s_{1}, 1_{S'}), \dots, (s_{\ell_2}, 
    1_{S'})$ is an irreducible $(\{0\}\times T')$-sequence of length $\ell_1 
    + \ell_2 = D(S', T') + D(S, \{0\}) - 2$, as desired.
  \end{proof}
\end{lemma}
\begin{definition}
  A semigroup $S$ is \textit{almost unit-stabilized} if $S$ has a 0 element 
  and $(S, \{0\})$ is a unit-stabilized pair.
\end{definition}
\begin{remark}
  If $S$ is unit-stabilized and contains a 0 element, the $S$ is also almost unit-stabilized.
\end{remark}
\begin{lemma}
  \label{auslem}
  Suppose $S$ is almost unit-stabilized, and suppose $S'$ is any other semigroup.
  Then
  \begin{equation*}
    D(S\times S') = \max(D(U(S)\times S'), D(S, \{0\}) + D(S') - 1).
  \end{equation*}
  \begin{proof}
    This follows directly from \Cref{relthm} and 
    \Cref{zerolemma}.
  \end{proof}
\end{lemma}
\begin{corollary}
  \label{auscor}
  Suppose $S$ is almost unit-stabilized and $\Delta(S) = 0$. Then for any
  semigroup $S'$, $D(S\times S') = D(U(S)\times S')$.
  \begin{proof}
    If $\Delta(S) = 0$, then $D(S, \{0\})\le D(U(S))$. Thus by \Cref{prodlemma}, 
    \begin{align*}
      D(S, \{0\}) + D(S') - 1 \le D(U(S)) + D(S') - 1\le D(U(S)\times S'),
    \end{align*}
    which means that by \Cref{auslem},
    \[
      D(S\times S') = \max(D(U(S) \times S'), D(S, \{0\}) + D(S') - 1) = D(U(S) \times S').
    \]
  \end{proof}
\end{corollary}
\begin{theorem}
  \label{austhm}
  Suppose $S_1, \dots, S_k$ are almost unit-stabilized semigroups. Suppose further that
  for $i > \ell$, $\Delta(S_i) = 0$. Then
  \begin{equation*}
    D(S_1\times\dots\times S_k) = \max_{I\subseteq [1,\ell]}\left[D\left(\prod_{i\in 
    [1,k]\backslash I} U(S_i)\right) + \sum_{i\in I} 
  (D(S_i, \{0\}) - 1)\right].
  \end{equation*}
  \begin{proof}
    First, by $\ell$ applications of \Cref{auslem}, we have that for any semigroup $S'$,
    \begin{equation*}
      D(S_1\times\dots\times S_\ell\times S')
      = \max_{I\subseteq [1,\ell]}\left[D\left(\left(\prod_{i\in
      [1,\ell]\backslash I} U(S_i)\right)\times S'\right) + \sum_{i\in I} (D(S_i,
    \{0\}) - 1)\right].
    \end{equation*}
    Now if we let $S' = S_{\ell + 1} \times \dots \times S_k$ and apply
    \Cref{auscor} $k - \ell$ times, we have
    \begin{equation*}
      D(S_1\times\dots\times S_k)
      = \max_{I\subseteq [1,\ell]}\left[D\left(\prod_{i\in [1,k]\backslash I}
      U(S_i)\right) + \sum_{i\in I} (D(S_i, \{0\}) - 1)\right].
    \end{equation*}
  \end{proof}
\end{theorem}
\begin{corollary}
  If $T = S_1 \times\dots\times S_k$ is a product of almost unit-stabilized 
  semigroups, then
  \begin{equation*}
    \Delta(T) \le \Delta(S_1)+\dots+\Delta(S_k).
  \end{equation*}
  \begin{proof}
    Applying \Cref{austhm} to the (trivial) case where $\ell = k$, we 
    have
    \begin{equation*}
      D(T) = \max_{I\subseteq [1,k]}\left[D\left(\prod_{i\in 
      [1,k]\backslash I} U(S_i)\right) + \sum_{i\in I} (D(S_i, \{0\}) - 
    1)\right].
    \end{equation*}
    However, using the inequality $D(G)\ge D(G / H) + D(H) - 1$ (a special case
    of \Cref{grplem} where $S$ is the trivial group), we have that for any
    $I\subseteq [1,k]$,
    \begin{align*}
      D\left(\prod_{i\in [1,k]\backslash I} U(S_i)\right) + \sum_{i\in I} 
      (D(U(S_i)) - 1) \le  D\left(\prod_{i\in [1,k]} U(S_i)\right) = 
      D(U(T)).
    \end{align*}
    Thus
    \begin{align*}
      D(T) &= \max_{I\subseteq [1,k]}\left[D\left(\prod_{i\in 
      [1,k]\backslash I} U(S_i)\right) + \sum_{i\in I} (D(S_i, \{0\}) - 
    1)\right] \\ &\le \max_{I\subseteq [1,k]}\left[D(U(T)) + \sum_{i\in I} 
    (D(S_i, \{0\}) - D(U(S_i))\right] \\ &= D(U(T)) + \sum_{i = 1}^k \max(0, 
    D(S_i, \{0\}) - D(U(S_i)))\\
    &= D(U(T)) + \sum_{i = 1}^k \Delta(S_i),
    \end{align*}
    from which the desired result immediately follows.
  \end{proof}
\end{corollary}
\begin{remark}
  It would be nice if it were true in general that $\Delta(S_1\times S_2) \le
  \Delta(S_1) + \Delta(S_2)$. However, this is not true, even when one of
  $S_1$ or $S_2$ is almost unit-stabilized! For an example, let $S_1 = \BZ /
  3\BZ$ and let $S_2 = S\times S'$, where $S' = \BZ/4\BZ$ and $S = \{0, 1, 2,
  4\}\subseteq \BZ/7\BZ$.

  We have that the following:
  \begin{itemize}
    \item $U(S_1)\simeq C_2, U(S)\simeq C_3, U(S')\simeq C_2$.
    \item $S_1, S, S'$ are unit-stabilized.
    \item $D(S_1) = D(U(S_1)) = 2$.
    \item $D(S) = D(U(S)) = 3$.
    \item $D(S') = D(S', \{0\}) = 3$. $D(U(S')) = 2$.
    \item $D(S_2) = D(S\times S') = D(U(S)\times S') = \max(D(U(S)\times
      U(S')), D(U(S)) + D(S', \{0\}) - 1) = \max(6, 3 + 3 - 1) = 6$,
      $D(U(S_2)) = D(C_6) = 6$.
  \end{itemize}
  However,
  \begin{align*}
    D(S_1\times S_2) = D(S_1\times S \times S') &= D(U(S_1)\times U(S)\times S') \\
    &= \max(D(C_6\times U(S')), D(C_6) + D(S', \{0\} - 1) \\
    &= \max(7, 6 + 3 - 1) = 8.
  \end{align*}
  On the other hand,
  \[
    D(U(S_1\times S_2)) = D(U(S_1)\times U(S)\times U(S')) = D(C_6\times C_2) = 7.
  \]
  Thus $\Delta(S_1) = \Delta(S_2) = 0$ but $\Delta(S_1\times S_2) = 1$. In
  fact, it is possible to construct $S_1, S_2$ such that $\Delta(S_1) =
  \Delta(S_2) = 0$ and one of $S_1, S_2$ is almost unit-stabilized, but
  $\Delta(S_1\times S_2)$ is arbitrarily large.
\end{remark}
\section{Local Rings}
In this section, we look at the case where $R$ is a (finite) local ring with
maximal ideal $\mf{m}$ and residue field $k$. It turns out local rings are all
either unit-stabilized or almost unit-stabilized. 
\begin{lemma}
  For all $a\in R$ there exists a positive integer $n$ such that $a\mf{m}^n = 0$.
  \begin{proof}
    Since $R$ is finite, the chain of ideals 
    \[
      aR\supseteq a\mf{m}\supseteq a\mf{m}^2\supseteq \cdots
    \]
    must stabilize, so $a\mf{m}^n = a\mf{m}^{n+1}$ for some positive integer $n$. Then by Nakayama's lemma we must have $a\mf{m}^n = 0$.
  \end{proof}
\end{lemma}
\begin{lemma}
  If $a\neq 0$ and $\Ann_R(a) = \Ann_R(ab)$, then $b\in U(R)$.
  \begin{proof}
    Suppose instead that $b\notin U(R)$, so $b\in \mf{m}$.
    Consider the minimum $n$ such that $a\mf{m}^n = 0$. (Clearly $n > 0$.) Then $ab\mf{m}^{n-1}
    \subseteq a\mf{m}^n = 0$ but $a\mf{m}^{n - 1}\neq 0$, which means that
    $\mf{m}^{n-1}\not\subseteq \Ann_R(a)$ but $\mf{m}^{n-1}\subseteq
    \Ann_R(ab)$, contradiction.
  \end{proof}
\end{lemma}
\begin{theorem}
  \label{goodrings}
  We have:
  \begin{enumerate}
    \item If $k \not \simeq \FT$, then $R$ is unit stabilized.
    \item If $k \simeq \FT$, then $R$ is almost unit-stabilized.
  \end{enumerate}
  \begin{proof}
    Note that for all $a\in R$, $\Stab_{U(R)}(a) = (1 + \Ann_R(a))
    \intersect U(R)$. In addition, if $a\neq 0$, then $\Ann_R{U(R)}$ is a
    proper ideal of $R$, which means that $\Ann_R{U(R)}\subseteq \mf{m}$, so
    $1 + \Ann_R{U(R)}\subseteq U(R)$, which means that $\Stab_{U(R)}(a) = 1
    + \Ann_R(a)$. Thus if $ab\neq 0$ and $\Stab_{U(R)}(a) =
    \Stab_{U(R)}(ab)$, then $\Ann_R(a) = \Ann_R(ab)$, which by the previous
    lemma implies  $b\in U(R)$.

    In the case where $R / \mf{m}\not\simeq \FT$, then we have $1 +
    \mf{m}\subsetneq U(R)$, so for all $a \neq 0$, $\Stab_{U(R)}(a)\subseteq
    1 + \mf{m} \subsetneq U(R) = \Stab_{U(R)}(0)$, which means that $R$ is
    unit-stabilized.
  \end{proof}
\end{theorem}
\begin{corollary}
    If $k\not \simeq \FT$, then $\Delta(R) = 0$.
\end{corollary}
\begin{theorem}
  \label{badrings}
  Suppose $R$ is a finite local ring. Then $\Delta(R) \le 1$, and equality holds if and only
  if $R$ is isomorphic to one of the following rings:
  \begin{itemize}
    \item $\BZ / 2\BZ$
    \item $\BZ / 4\BZ$
    \item $\BZ / 8\BZ$
    \item $\FT[x] / (x^2)$
  \end{itemize}
  \begin{proof}
    We have already taken care of the case where $k\not\simeq\FT$.
    Suppose $R$ is a local ring with residue field $\FT$. Then $\abs{R} = 2^n$ for a positive integer $n$.
    We will show the following: 
    \begin{enumerate}
      \item $D(R, 0)\le n + 1$.
      \item If $D(R, 0) = n + 1$, then $\mf{m}^{n - 1}\neq 0$.
      \item $D(U(R))\ge n$.
      \item If $D(U(R)) = n$, then $U(R)\simeq C_2^{n-1}$.
      \item If $n \ge 3$ and $\mf{m}^{n - 1}\neq 0$ and $U(R)\simeq C_2^{n
        - 1}$, then $\mf{m} = (2)$ and $R\simeq \BZ / 2^n\BZ$.
      \item If $n\ge 4$, then $D(U(\BZ / 2^n\BZ)) \not \simeq C_2^{n -
        1}$.
      \item If $n\le 2$, then $R\simeq \BZ / 2\BZ, \BZ / 4\BZ$, or $\FT[x]
        / (x^2)$.
    \end{enumerate}
    To finish the proof, note that since $R$ is almost unit-stabilized,
    $D(R) = \max(D(U(R)), D(R, 0)$, so from (1) and (2) we immediately get
    $\Delta(R)\le 1$. Now suppose $\Delta(R) = 1$. Since $D(R, 0) \le n + 1$,
    $D(U(R))\ge n$, and $\Delta(R) = \max(0, D(R, 0) - D(U(R))$, equality
    must hold in both cases. Thus by (2), $\mf{m}^{n - 1}\neq 0$; by (4),
    $U(R)\simeq C_2^{n - 1}$; by (5), this means that either $\abs{R}\le 4$
    or $R\simeq \BZ / 2^n\BZ$. In the first case, by (7), $R$ must be either
    $\BZ / 2\BZ, \BZ / 4\BZ$, or $\FT[x] / (x^2)$. Otherwise, by (6), we
    must have $n = 3$, in which case $R\simeq \BZ / 8\BZ$. 

    Finally, we have the following values of $D(R)$ and $D(U(R))$ for $R = \BZ / 2\BZ, \BZ /
    4\BZ, \BZ / 8\BZ, \FT[x] / (x^2)$
    \begin{center}
      \begin{tabular}{c|c|c}
        $R$ & $D(R)$ & $D(U(R))$ \\
        \hline
        $\BZ / 2\BZ$ & 2 & 1 \\
        $\BZ / 4\BZ$ & 3 & 2 \\
        $\BZ / 8\BZ$ & 4 & 3 \\
        $\FT[x] / (x^2) $ & 3 & 2 \\
      \end{tabular}
    \end{center}
    In each of these cases, we have $\Delta(R) = 1$. 

    Here are proofs of the 7 claims:
    \begin{enumerate}
      \item Suppose $a_1, \dots, a_{n+1}$ is an irreducible sequence in
        $R$ with $a_1 \cdots a_{n + 1} = 0$. Then none of $a_1, \dots,
        a_{n+1}$ are in $U(R)$, or else they could be omitted. In
        addition, $a_1\cdots a_n\neq 0$. Thus if $p_i = a_1\cdots a_i$, we
        have $\Stab_{U(R)}(p_i)\subsetneq \Stab_{U(R)}(p_{i+1})$ for
        $0\le i< n$. In particular, this means that
        $\abs{\Stab_{U(R)}(p_i)}\ge 2^i$. However, when $i = n$ this
        gives $2^{i - 1} = \abs{U(R)}\ge \abs{\Stab_{U(R)}(p_i)} \ge
        2^i$, contradiction.
      \item Suppose $D(R, 0) = n + 1$, and let $a_1, \dots, a_n$ be an
        irreducible sequence in $R$ with $a_1\cdots a_n = 0$. Then none
        of the $a_i$ are units, so $a_i\in \mf{m}$. In addition,
        $a_1\cdots a_{n-1} \neq 0$, which means that $\mf{m}^{n - 1}\neq
        0$.
      \item We have that $\abs{U(R)} = 2^{n - 1}$. Thus we can write 
        \[
          U(R) = C_{2^{e_1}}\times \cdots \times C_{2^{e_r}},
        \]
        where $e_1 +
        \dots + e_r = n - 1$. Since $R$ is a 2-group, by
        \cite{Olson1969}, 
        \[
          D(U(R)) = 1 + \sum_{i = 1}^r (D(C_{2^{e_i}}) - 1) = 1 +
          \sum_{i = 1}^r (2^{e_i} - 1)
          \ge 1 + \sum_{i = 1}^r e_i = n.
        \]
        Here we also used the inequality $2^n - 1\ge n$ for all $n\in \BZ$.
      \item In the above equation, equality holds if and only if each of the $e_i$ is
        equal to 1, in which case $U(R)\simeq C_2^{n - 1}$.
      \item 
        Since $\mf{m}^{n- 1} \neq 0$, by Nakayama's lemma we have that $\mf{m}^i\subsetneq
        \mf{m}^{i - 1}$ for $1\le i\le n$, in which case
        $2\abs{\mf{m}^i}\le \abs{\mf{m}^{i - 1}}$. Thus we have
        \[
          2^{n} = \abs{\mf{m}^0} \ge 2\abs{\mf{m}^1} \ge \dots \ge
          2^{n - 1} \abs{\mf{m}^{n - 1}} \ge 2^{n -1 } \cdot 2 = 2^n,
        \]
        so equality holds in each of the above inequalities. In
        particular, we have $\abs{\mf{m} / \mf{m^2}} = 2$, so
        \[
          \dim_k(\mf{m} / \mf{m^2}) = 1,
        \]
        which means that $\mf{m}$ is principal.

        Suppose $\mf{m} = (t)$. Then $t - 1$ is a unit, and since
        $U(R)\simeq C_2^{n - 1}$, we must have $(t - 1)^2 = 1$, or $t^2
        = 2t$.
        
        I claim that $2\in \mf{m} \backslash \mf{m}^2$. Suppose otherwise, that $2\in \mf{m}^2$. Then we have $t^2 = 2t\in \mf{m}^3$,
        which means that $\mf{m}^2 = (t^2) = \mf{m}^3$. By Nakayama's
        lemma, this means that $\mf{m}^2 = 0$, which also means $\mf{m}^{n - 1} = 0$ as $n\ge 3$, 
        contradiction.

        Thus $2\in \mf{m} \backslash \mf{m}^2$, which means that $2
        = ut$ for some unit $u$, so $\mf{m} = (t) = (2)$. Finally,
        note that 
        \[
          (2^{n - 1}) = \mf{m}^{n - 1} \neq 0,
        \] so $2^{n -
        1}\neq 0$. Thus in the additive group of $R$, 1 has order $2^n$, which 
        means that we must have $R\simeq
        \BZ / 2^n\BZ$.
      \item If $n\ge 4$, then 5 is a unit in $\BZ / 2^n\BZ$, but $3^2\not
        \equiv 1\pmod{2^n}$, so $U(\BZ / 2^n\BZ)\not\simeq C_2^{n - 1}$.
      \item When $n = 1$, there is a unique ring with two elements, $\BZ /
        2\BZ$. Now if $n = 2$, $\abs{\mf{m}} = 2$, so $\mf{m} = \{0,
        x\}$ for some $x\neq 0, 1$. Then $R = \{0, 1, x, x - 1\}$. We
        have two cases:
        \begin{itemize}
          \item $2\neq 0$. Then we must have $x = 2$ and $x - 1 =
            3$, so $R\simeq \BZ / 4\BZ$.
          \item 2 = 0. Then $1 = (x + 1)^2 = x^2 + 1$, so $x^2 = 0,
            (x+1)^2 = 1, x(x+1) = x$. In this case, $R\simeq \FT[x]
            / x^2$.
      \end{itemize}
    \end{enumerate}
  \end{proof}
\end{theorem}
\section{Proof of \texorpdfstring{\Cref{thm1}}{Theorem \ref{thm1}}}
Let $R_1 = \BZ / 2\BZ, R_2 = \BZ / 4\BZ, R_3 = \BZ / 8\BZ, R_4 = \FT[x]/(x^2)$ denote the ``bad'' local rings.
We have that $R$ is of the form
\[
  R_1^{k_1} \times R_2^{k_2} \times R_3^{k_3} \times R_4^{k_4
  } \times R',
\]
where $R'$ is a product of local rings not isomorphic to $R_1, R_2, R_3$, or $R_4$.
Then by \Cref{badrings}, $R'$ is a product of almost unit-stabilized local rings with $\Delta = 0$. In addition, we have
\begin{itemize}
  \item $U(R_1) = 0$, $D(R_1) = 1$.
  \item $U(R_2)\simeq C_2$, $D(R_2) = 2$.
  \item $U(R_3)\simeq C_2 \times C_2$, $D(R_3) = 3$.
  \item $U(R_4)\simeq C_2$, $D(R_4) = 2$.
\end{itemize}
Plugging this all into \Cref{austhm}, we have
\begin{align*}
  D(R) 
  &= \max_{n_i\in [0, k_i]} \bigg[
    D\left( {C_2}^{(k_2 - n_2) + 2(k_3 - n_3) + (k_4 - n_4)} \times U(R')\right) + n_1 + 2n_2 + 3n_3 + 2n_4
  \bigg] \\
  &= \max_{n_i\in [0, k_i]} \bigg[
    D\left( {C_2}^{k_2 + k_4 + 2k_3 - (n_2 + n_4) - 2n_3} \times U(R')\right) + k_1 + 2(n_2 + n_4) + 3n_3
  \bigg] \\
  &= \max_{0\le a\le k_2 + k_4, 0\le b\le k_3} \bigg[
    D\left( U(R') \times {C_2}^{k_2 + k_4 + 2k_3 - a - 2b} \right) + k_1 + 2a + 3b
  \bigg].
\end{align*}
\section{A few more corollaries of \texorpdfstring{\Cref{thm1}}{Theorem \ref{thm1}}}
\begin{theorem}
  \label{randomthms}
  $D(R) = D(U(R))$ in any of the following scenarios:
  \begin{enumerate}
      \item $\abs{R}$ is odd.
      \item $R$ is an $\BF_q$ algebra, where $q$ is a power of 2 other than 2.
      \item $R = \BZ / n\BZ$, where $16\mid n$.
      \item $R = \FT[x] / (f)$, where $v_x(f), v_{x + 1}(f)\notin \{1, 2\}$.
      \item $R$ is of the form $\BZ / 4\BZ \times R'$ or $\FT[x] / (x^2) \times R'$, where $\abs{U(R')} > 1$ is odd and $\Delta(R') = 0$.
      \item $R = \FT[x] / (x^2g)$, where $g$ is a product of (at least one) distinct irreducibles that are not $x$ or $x + 1$.
  \end{enumerate}
  \begin{proof}
    \ 
    \begin{enumerate}
      \item This follows immediately from \Cref{thm3}, as a ring with odd order cannot have any ideals of index 2.
      \item Note that any quotient ring of $R$ is also an $\BF_q$ algebra, which means that any proper ideal has index at least $q$ in $R$. Thus $n_2(R) = 0$, so by \Cref{thm3}, $\Delta(R) = 0$.
      \item In the setting of \Cref{thm2}, it is easy to see that for the ring $R = \BZ / n\BZ$ with $16\mid n$, so we have $k_1 = k_2 = k_3 = k_4 = 0$, which means that $\Delta(R) = 0$ as well.
      \item This is a direct consequence of \Cref{polythm} below, which is in turn a consequence of \Cref{thm2}.
      \item We first need the following lemma:
        \begin{lemma}
          Suppose $G$ is a nontrivial group of odd order. Then $D(G\times C_2) > D(G) + 1$.
          \begin{proof}
            (We will use additive notation for this proof).
            Let $(g_1, \dots, g_\ell)$ be an irreducible
            sequence in $G$, where $\ell = D(G) - 1$. Note that the map $g\mapsto g + g$ in $G$ is injective, which means that the map $(\times \frac12)$ is well-defined. Then the sequence
            \[
              (g_1, 0), \dots, (g_{\ell - 1}, 0),
              (\frac12 g_\ell, 1),
              (\frac12 g_\ell, 1),
              (\frac12 g_\ell, 1)
            \]
            is an irreducible sequence in $G\times C_2$ of length $\ell + 2 = D(G) + 1$, so $D(G\times C_2) > D(G) + 1$.
          \end{proof}
        \end{lemma}
        Since $\abs{U(R')}$ is odd, the decomposition of $R'$ into a product of local rings cannot have any rings isomorphic to $\BZ / 4\BZ, \BZ / 8\BZ$, or $\FT[x] / (x^2)$. In addition, since $\Delta(R') = 0$, by \Cref{thm2}, this decomposition cannot have any rings isomorphic to $\BZ / 2\BZ$. Thus the setting of \Cref{thm1} applies; plugging in both $R = \BZ / 4\BZ \times R'$ and $\FT[x] / (x^2) \times R'$ gives
        \[
          D(R) = \max(D(U(R')) + 2, D(U(R') \times C_2))
        \]
        However, by the lemma, $D(U(R') \times C_2) > D(U(R')) + 1$, which means that $D(R) = D(U(R') \times C_2) = D(U(R))$.
      \item 
        If $f\in \FT$ is irreducible, then the unit group of $\FT[x] / (f)$ has odd order (in particular, it has order $2^{\deg f} - 1$.) Thus if $g$ is a product of irreducibles other than $x, x + 1$, then $\abs{U(\FT[x] / (g))} > 1$ is odd, and by \Cref{thm2}, $\Delta(\FT[x] / (g)) = 0$. Thus this is a special case of (5), with $R' = \FT[x] / (g)$ and $R\simeq \FT[x] / (x^2) \times R'$, so $\Delta(R) = 0$.
    \end{enumerate}
  \end{proof}
\end{theorem}
Note that (1) in \Cref{randomthms} is a generalization of \Cref{wangguo} and (2) is a generalization of \Cref{wang}. In addition, we have the following refinement of \Cref{qwz}.
\begin{theorem}
  \label{polythm}
  Let $f\in \FT[x]$ be nonconstant and $R = \FT[x] / (f)$. If $f = x^a (x +
  1)^b g$, where $\gcd(g, x(x+1)) = 1$, then 
  \[
    \delta_{a1} + \delta_{b1}\le D(R) - D(U(R)) \le \delta_{a1} + \delta_{b1} +
    \delta_{a2} + \delta_{b2},
  \]
  where $\delta_{ij}$ is the Kronecker $\delta$ function.
  \begin{proof}
    Note that if $f$ factors into irreducibles as $f = f_1^{e_1}\cdots f_n^{e_n}$, then the decomposition of $\FT[x] / (f)$ into products of local rings is
    \[
      \FT[x]/(f) \simeq \FT[x] / (f_1^{e_1}) \times \dots \times \FT[x] / (f_n^{e_n}).
    \]
    In the setting of \Cref{thm2}, we have $k_1 = \delta_{a1} + \delta_{b1}$, $k_2 = k_3 = 0$, and $k_4 = \delta_{a2} + \delta_{b2}$. The desired result immediately follows.
  \end{proof}
\end{theorem}
Finally, we give a partial answer to the general problem stated by Wang and Guo in \cite{WangGuo2008}:
\begin{theorem}
  If $R = \BZ/n_1\BZ \times \dots \times \BZ/n_r\BZ$, then
  \[
    \#\{1\le i\le r: 2\| n_i\} \le D(R) - D(U(R)) \le \#\{1\le i\le r: 2\mid n_i, 16\nmid n_i\}.
  \]
  In addition, equality holds on the right when all of the $n_i$ are powers of 2.
  \begin{proof}
    Again in the setting of \Cref{thm2}, we have $k_4 = 0$, and for $1\le i\le 3$, 
    \[
      k_i = \#\{1\le i\le r: 2^i\| n_i\}
    \]
    Thus by \Cref{thm2},
    \[
      D(R) - D(U(R)) \ge k_1 = \#\{1\le i\le r: 2\| n_i\}
    \]
    and 
    \begin{align*}
      D(R) - D(U(R))&\le k_1 + k_2 + k_3 + k_4 \\
      &= \sum_{i = 1}^3 \#\{1\le i\le r: 2^i\| n_i\}
      \\
      &= \#\{1\le i\le r: 2\mid n_i, 16\nmid n_i\}.
    \end{align*}
    In addition, if each of the $n_i$ is a power of 2, $U(R)$ is a 2-group, so by \Cref{thm2}, equality holds on the right.
  \end{proof}
\end{theorem}
\section{Further Direction}
Unfortunately, we still do not have a complete classification of the finite rings $R$ for which $D(R) = D(U(R))$, even if we restrict to the simple cases where $R$ is of the form $\BZ / n\BZ$ or $\FT[x] / (f)$.

As we can see, these questions depend on the relation between $D(G)$ and $D(G\times C_2^e)$ for an abelian group $G$. This is what we know:

If $G = C_{n_1} \times \dots \times C_{n_r}$, where $n_1 \mid \cdots \mid n_r$, then we can define $M(G) = 1 + \sum_{i = 1}^r (n_r - 1)$. Clearly $D(G)\ge M(G)$. However, we have the following:
\begin{proposition}
  \label{gooddelta}
$D(G) = M(G)$ in each of the following cases:
\begin{enumerate}
  \item $G$ is a $p$-group. \cite{Olson1969}
  \item $G = C_n \times C_{kn}$. \cite{Olson1969}
  \item $G = C_{p^ak} \times H$, where $H$ is $p$-group and $p^a \ge M(H)$. \cite{vebii}
  \item $G = C_2 \times C_{2n} \times C_{2kn}$, where $n, k$ are odd and the largest prime dividing $n$ is less than 11. \cite{vebii}
  \item $G = C_2^3 \times C_{2n}$, where $n$ is odd. \cite{baayen}
\end{enumerate}
\end{proposition}
(For a more complete list, see \cite{mazur}.) On the other hand, it is not true that $D(G) = M(G)$ for all $G$.
\begin{proposition}
  \label{baddelta}
  $D(G) > M(G)$ if $G = C_2^e \times C_{2n} \times C_{2nk}$ for $n, k$ odd and $e\ge 3$. \cite{Geroldinger}
\end{proposition}
In what follows, we will try and compute exact values of $\Delta(R)$ when $R$ is of the form $\BZ / n\BZ$ or $\FT[x] / (f)$.
\subsection{\texorpdfstring{$R = \BZ / n\BZ$}{R = Z / nZ}}
By \Cref{thm3}, we have $\Delta(\BZ / n\BZ) \le 1$.
\begin{question}
  For $n\in\BZ^+$, when is $\Delta(\BZ / n\BZ)$ equal to 0 and when is it 1?
\end{question}
Note that we have already resolved this problem for most $n$. In particular, we have the following:
\begin{itemize}
  \item If $n$ is odd or $16 \mid n$, $\Delta(R) = 0$.
  \item If $n = 2b$ where $b$ is odd, then $\Delta(R) = 1$.
  \item If $n = 4b$ where $b$ is odd, then by \Cref{thm1} we have
    \[
      \Delta(R) = \max(0, D(U(\BZ / b\BZ)) + 2 - D(U(\BZ / n\BZ)))
    \]
    However, $U(\BZ / n\BZ) = U(\BZ / b\BZ) \times C_2$, which means that 
    \[
      D(U(\BZ / n\BZ)) \ge D(U(\BZ / b\BZ)) + 1.
    \]
    Thus $\Delta( \BZ / n\BZ) = 1$ if and only if 
    \[
      D(U(\BZ / n\BZ)) = D(U(\BZ / b\BZ)) + 1.
    \]
    Using this, we can calculate $\Delta(R)$ in the following special cases using \Cref{gooddelta,baddelta} in conjunction with \Cref{thm1}:
    \begin{itemize}
      \item If $b$ is a prime power then $\Delta(R) = 1$.
      \item If $b$ is a product of distinct primes of the form $2^{2^n} + 1$, then $\Delta(R) = 1$.
      \item If $b = p_1^{e_1} p_2^{e_2}$ and $\gcd(\phi(p_1^{e_1}), \phi(p_2^{e_2}))$ has no prime factors greater than 7, then $\Delta(R) = 1$.
      \item If $b = p_1^{e_1}p_2^{e_2}p_3^{e_3}$ and $
        \phi(p_1^{e_1})/2,
        \phi(p_2^{e_2})/2,
        \phi(p_3^{e_3})/2
        $ are pairwise relatively prime, then $\Delta(R) = 1$.
      \item If $b = p_1^{e_1}p_2^{e_2}p_3^{e_3}p_4^{e_4}$, $p_1, p_2, p_3, p_4 \equiv 3\pmod{4}$, and $
        \phi(p_1^{e_1})/2,
        \phi(p_2^{e_2})/2,
        \phi(p_3^{e_3})/2,
        \phi(p_4^{e_4})/2
        $ are pairwise relatively prime, then $\Delta(R) = 0$.
    \end{itemize}
  \item If $n = 8b$ where $b$ is odd, then using similar logic as the previous case we have $\Delta(\BZ / n\BZ) = 1$ if and only if \[
      D(U(\BZ / n\BZ)) = D(U(\BZ / b\BZ)) + 2.
    \]
    Thus we can calculate $\Delta(R)$ in the following special cases:
    \begin{itemize}
      \item If $b$ is a prime power then $\Delta(R) = 1$.
      \item If $b$ is a product of distinct primes of the form $2^{2^n} + 1$, then $\Delta(R) = 1$.
      \item If $b = p_1^{e_1} p_2^{e_2}$ and $\gcd(\phi(p_1^{e_1}), \phi(p_2^{e_2})) = 1$, then $\Delta(R) = 1$.
      \item If $b = p_1^{e_1}p_2^{e_2}p_3^{e_3}$, $p_1, p_2, p_3\equiv 3\pmod{4}$, and 
        \[
          \gcd[
          \phi(p_1^{e_1})/2,
          \phi(p_2^{e_2})/2,
          \phi(p_3^{e_3})/2] = 2
        \]
        and for $i\neq j$, 
        $ \gcd[ \phi(p_i^{e_i})/2, \phi(p_j^{e_j})/2] $ has
        no prime factors greater than 7, then $\Delta(R) = 0$.
      \item If $b = p_1^{e_1}p_2^{e_2}p_3^{e_3}p_4^{e_4}$, $p_1, p_2, p_3, p_4 \equiv 3\pmod{4}$, and $
        \phi(p_1^{e_1})/2,
        \phi(p_2^{e_2})/2,
        \phi(p_3^{e_3})/2,
        \phi(p_4^{e_4})/2
        $ are pairwise relatively prime, then $\Delta(R) = 0$.
    \end{itemize}
\end{itemize}
However, the general cases $n = 4b, 8b$ remain open.
\subsection{\texorpdfstring{$R = \FT[x]/(f)$}{R = F2[x]/(f)}}
In this case, \Cref{thm3} gives the bound $\Delta(\FT[x] / (f)) \le 2$.
\begin{question}
  For $f\in \FT[x]$, when is $\Delta(\FT[x] / (f))$ equal to 0, when is it 1, and when is it 2?
\end{question}
Again, we have already answered this question for most $f\in \FT[x]$. Let $f = x^a(x+1)^b g$, where $\gcd(g, x(x + 1)) = 1$. Then
\begin{itemize}
  \item If $g = 1$, then $\Delta(R) = \delta_{a1} + \delta_{a2} + \delta_{b1} + \delta_{b2}$.
  \item If $a, b\neq 2$, then $\Delta(R) = \delta_{a1} + \delta_{b1}$.
  \item If $a = 2$ and $b \le 1$ and $\abs{U(\FT[x] / (g))}$ is odd, then $\Delta(R) = b$.
    (Note that this condition corresponds to $g$ being a product of distinct irreducibles.)
  \item If $a = b = 2$ and $U(\FT[x] / (g))$ is cyclic (i.e. product of distinct irreducibles with pairwise relatively prime degree), then $\Delta(R) = 1$.
  \item If $a = 2$, $b\in [3, 14]\cup [17, 24] \cup[33, 40]$ and $U(\FT[x]/(g))$ is cyclic, then $\Delta(R) = 0$.
\end{itemize}
\subsection{Recap}
For general $n\in \BZ$, the rank of $U(\BZ / n\BZ)$ is equal to the number of prime factors of $n$, which gets arbitrarily large. 

Unfortunately, very little is known about $D(G) - M(G)$ for general $G$ of large rank. 

In a similar vein, the rank of $\FT[x] / (f)$ for general $f\in \FT[x] / (f)$ can get annoyingly large, and even the special case $\Delta(\FT[x]/(x^2g))$ with $g\in \FT[x]$ irreducible is tricky to calculate.

For example, \cite{Geroldinger} gives for any $k > 0$, $D(C_{2k + 1} \times {C_2}^r) =
4k + 1 + r$ for $1\le r\le 4$ but $D(C_{2k + 1} \times {C_2}^5) > 4k + 7$.

If $f$ is an irreducible polynomial in $\FT[x]$ of degree $d$, one
can check that $U(\FT[x] / (f^2))\simeq C_{2^d - 1}\times {C_2}^d$.
Thus if $R_d = \FT[x] / (x^2f_d^2)$, where $f_d$ is an irreducible polynomial of
degree $d$, applying \Cref{thm1} gives 
\[
  D(R_d) = \max(D(C_{2^d - 1} \times {C_2}^d) + 2, D(C_{2^d - 1} \times
  {C_2}^{d+1}))
\]
which gives
$\Delta(R_2) = \Delta(R_3) = 1$ but $\Delta(R_4) = 0$. The general problem of determining $\Delta(\FT[x] / (f))$ for non-squarefree $f\in \FT[x]$ is still open for reasons such as this. For more information about the growth of $D(C_{2k + 1} \times {C_2}^d)$, we refer the reader to \cite{mazur}.
\subsection{Infinite Semigroups}
Finally, we look at the case where $S$ is an infinite commutative semigroup.
By \Cref{thm4} and \Cref{inflemma}, when $S$ is either
\begin{itemize}
  \item a group, or
  \item the semigroup of a commutative ring $R$ under multiplication,
\end{itemize}
then $D(S)$ being finite implies $S$ is finite. However, there are semigroups $S$ such that $S$ is infinite but $S$ is finite.
Thus we have the following question:
\begin{question}
  For what other families $\mc{F}$ of (commutative, unital) semigroups is it true $S\in \mc{F}$ and $\abs{S} = \infty$ implies $D(S) = \infty$?
\end{question}
\section*{Acknowledgments}
This research was conducted at the University of Minnesota Duluth REU and was 
supported by NSF grant 1358695 and NSA grant H98230-13-1-0273. The author thanks 
Joe Gallian for suggesting the problem and supervising the research, and Benjamin Gunby for 
helpful comments on the manuscript.

\bibliographystyle{plain}
\bibliography{writeup}
\end{document}